\documentclass[11pt]{article}

\usepackage{IEEEtrantools}
\usepackage{amsfonts}
\usepackage{amssymb}
\usepackage{amsmath,amsfonts,amsthm}
\usepackage{mathrsfs}
\usepackage{url}
\usepackage{tikz}
\usetikzlibrary{matrix,arrows,decorations.pathmorphing}

\usepackage{color}

\newcommand{\bbbt}{\mathbb{T}}

\newcommand{\be}{\begin{equation}}
\newcommand{\ee}{\end{equation}}
\newcommand{\bea}{\begin{eqnarray}}
\newcommand{\eea}{\end{eqnarray}}
\newcommand{\bean}{\begin{eqnarray*}}
\newcommand{\eean}{\end{eqnarray*}}
\newcommand{\brray}{\begin{array}}
\newcommand{\erray}{\end{array}}
\newcommand{\biearray}{\begin{IEEEarray}{rCl}}
\newcommand{\eiearray}{\end{IEEEarray}}

\newcommand{\newsection}[1]{\setcounter{equation}{0}
\setcounter{dfn}{0}
\section{#1}}

\newcommand{\newsubsection}[1]{
\subsection{#1}}

\newtheorem{dfn}{Definition}[section]
\newtheorem{thm}[dfn]{Theorem}
\newtheorem{lmma}[dfn]{Lemma}
\newtheorem{ppsn}[dfn]{Proposition}
\newtheorem{crlre}[dfn]{Corollary}
\newtheorem{xmpl}[dfn]{Example}
\newtheorem{rmrk}[dfn]{Remark}

\newcommand{\bdfn}{\begin{dfn}\rm}
\newcommand{\bthm}{\begin{thm}}
\newcommand{\blmma}{\begin{lmma}}
\newcommand{\bppsn}{\begin{ppsn}}
\newcommand{\bcrlre}{\begin{crlre}}
\newcommand{\bxmpl}{\begin{xmpl}}
\newcommand{\brmrk}{\begin{rmrk}\rm}

\newcommand{\edfn}{\end{dfn}}
\newcommand{\ethm}{\end{thm}}
\newcommand{\elmma}{\end{lmma}}
\newcommand{\eppsn}{\end{ppsn}}
\newcommand{\ecrlre}{\end{crlre}}
\newcommand{\exmpl}{\end{xmpl}}
\newcommand{\ermrk}{\end{rmrk}}

\newcommand{\bbc}{\mathbb{C}}

\newcommand{\bbn}{\mathbb{N}}
\newcommand{\bbr}{\mathbb{R}}

\newcommand{\prf}{\noindent{\it Proof\/}: }

\def \qed { \mbox{}\hfill
$\Box$\vspace{1ex}}

\begin{document}


 \author{\sc  {Partha Sarathi Chakraborty, Bipul Saurabh}}
 \title{Gelfand-Kirillov dimension of some simple unitarizable  modules}
 \maketitle
 \begin{abstract} 
 Let  $\mathcal{O}_q(G)$ be the quantized algebra of regular functions on  a semisimple simply connected  compact Lie group $G$.
 Simple unitarizable left $\mathcal{O}_q(G)$-module are classified. 
 In this article  we compute  their Gelfand-Kirillov dimension where $G$ is of the type $A$, $C$ and $D$.
\end{abstract}

{\bf AMS Subject Classification No.:} {\large 16}P{\large 90}, {\large
17}B{\large 37}, {\large
  20}G{\large 42}

\newsection{Introduction}   
Earlier in sixties,  Gelfand and Kirillov (\cite{GelKir-1966aa}) introduced a notion 
of growth of an algebra or a module, namely Gelfand Kirillov dimension.   This invariant played a crucial role in the classification of 
Weyl algebras. It has been computed on many occasions (see  \cite{BueJimPas-1997aa}, \cite{McC-1989aa},  \cite{Tor-1999aa}, \cite{Jos-1978aa}). 
The Gelfand-Kirillov dimension of the universal enveloping algebra of a
finite dimensional Lie algebra is same as its dimension (Proposition 8.1.15,(iii) in \cite{McCRob-2001aa}). Even though this has
also been computed for some quantized universal enveloping algebras (\cite{McC-1989aa}) it is generally agreed
upon that it is not easy to calculate them for modules in general  (\cite{BueJimPas-1997aa}). Instead of considering the
quantized universal enveloping algebras one can also consider their duals namely the quantized algebra of regular functions  $\mathcal{O}_q(G)$ and their modules.  However in this case most of the literature is concerned with the type A situation only (example (4.8) page 258 in \cite{BueTorVer-2003aa}). Here we take up the cases of type C and D as well.

More specifically  given  a semisimple simply connected compact Lie group $G$, Korogodski and Soibelman (\cite{KorSoi-1998aa}) 
classified  all simple unitarizable modules   of the quantized algebra of regular functions on $G$ denoted by $\mathcal{O}_q(G)$. Given an element $w$ in the Weyl group and 
an element in a fixed maximal torus of $G$, they produced a  simple unitarizable $\mathcal{O}_q(G)$-module $V_{t,w}$ and showed that in each equivalence class there is exactly one such module. Some natural questions arise;
what is Gelfand Kirillov dimension of $\mathcal{O}_q(G)$-module $V_{t,w}$  and how is it related to the parameters $t$ and $w$. 
we prove that  Gelfand Kirillov dimension of $V_{t,w}$ is equal to  length of the element $w$  of  the Weyl group  provided $G$ is 
of type $A$, $C$ or $D$.  The key idea is to decompose the Weyl word $w$ into smaller parts in a certain way, prove the result for the last part and then use backward induction. This method also applies to the modules of quantized algebra of regular functions on certain  homogeneous spaces.

The paper is organized as follows. Next section describes all simple unitarizable modules of Hopf $*$-algebra $\mathcal{O}_q(G)$ and associates a diagram 
to each such module. Also, a brief introduction of the Weyl groups  of type $A$, $C$ and $D$ and a reduced expression of an element of the Weyl group are  given.
In the third section, we prove our main result. 
In the final section, we take certain  homogeneous spaces and prove similar result.

  Throughout the paper algebras are assumed to be unital and  over the field $\bbc$.
Elements of a Weyl group 
will be called Weyl words. We denote by $\ell(w)$  
the length of the Weyl word $w$.  
Let us denote 
by $\left\{e_n: n\in \bbn \right\}$ the standard basis of the 
vector space $c_{00}(\bbn)$ where  $c_{00}(\bbn)$ is  the set of finitely supported sequence of complex numbers. 
The endomorphism $e_n \mapsto e_{n-1}$ of $c_{00}(\bbn)$ will be denoted by   $S$.  The map $e_n \mapsto ne_n$ will be denoted
by  $N$. Given two  endomorphisms $T$ and $T^{'}$   of $c_{00}(\bbn)^{\otimes k}$  and  a subspace $V$  of 
$c_{00}(\bbn)^{\otimes k}$  we say that 
$T \sim T^{'}$ on $V$ if there exist integers $m_1,m_2, \cdots ,m_{k}$ and a nonzero constant $C$ such that 
\[
 T=CT^{'} (q^{m_1N}\otimes q^{m_2N}\otimes \cdots \otimes q^{m_{k}N})
\]
on $V$. Throughout this paper, $q$ will denote a real number in the interval  $(0,1)$ and $C$ is used to denote a generic constant.

\newsection{Simple unitarizable modules of quantized function algebras}
In this section, we recall the definition of quantized algebra of regular functions on a simply connected semisimple compact Lie group $G$ and the classification of its simple unitarizable modules. 
For a detailed treatment, we refer the reader to (\cite{KliSch-1997aa}, Chapter $3$ in \cite{KorSoi-1998aa}). 
Let $G$ be a simply connected semisimple compact Lie group, $\mathfrak{g}$ its complexified Lie algebra of rank $n$. Fix a nondegenerate symmetric ad-invariant  form $\left\langle\cdot,\cdot\right\rangle$ on $\mathfrak{g}$
such that its restriction to the real Lie algebra of $G$ is negative definite.
Let $\Pi:=\{ \alpha_1,\alpha_2, \cdots, \alpha_n \}$   be the set of simple roots.                 
For simplicity, we write the root $\alpha_i$ as $i$ and the reflection $s_{\alpha_i}$ defined by the root $\alpha_i$ as $s_i$. 
The  Weyl group $W_n$ of $G$ can be described as the group generated 
 by the reflections $\{s_{i}: 1\leq i \leq n\}$. 
\bdfn  \label{chap2-d-G_q}
Let $U_q(\mathfrak{g})$ be the quantized universal enveloping algebra  equipped with the
$*$-structure corresponding to the compact real form of $\mathfrak{g}$ (see page $161$ and $179$, \cite{KliSch-1997aa}).  
Then the Hopf $*$-subalgebra of the dual Hopf $*$-algebra of  $U_q(\mathfrak{g})$ consisting of matrix co-efficients
of finite dimensional unitarizable $U_q(\mathfrak{g})$-modules is called 
the quantized algebra of regular functions on $G$ (see page $96-97$, \cite{KorSoi-1998aa}).   We denote it by $\mathcal{O}_q(G)$. 
\edfn
Let $(\!(u_{j,\mathfrak{g}}^i)\!)$ be the defining corepresentation of $\mathcal{O}_q(G)$ if $G$ is of type $A_n$ and $C_n$ and  
the irreducible corepresentation of $\mathcal{O}_q(G)$ corresponding to the highest weight $(1,0,0,\cdots ,0)$ if $G$ is of type  $D_n$. 
In first case, entries of the matrix $(\!(u_{j,\mathfrak{g}}^i)\!)$ generate the Hopf $*$-algebra $\mathcal{O}_q(G)$. 
In latter case, they genarate a proper Hopf $*$-subalgebra of $\mathcal{O}_q(\mbox{Spin}(2n))$ which we denote as $\mathcal{O}_q(SO(2n))$.  
The generators of  $\mathcal{O}_q(\mbox{Spin}(2n))$ are the  matrix entries of the corepresentation 
$(\!(z_j^i)\!)$ of $\mathcal{O}_q(\mbox{Spin}(2n))$  with highest weight $(1/2,1/2,\cdots ,1/2)$. 
We denote the dimension of the corepresentation $(\!(u_{j,\mathfrak{g}}^i)\!)$ by $N_n$.  We will 
drop the subscript $\mathfrak{g}$ in $(\!(u_{j,\mathfrak{g}}^i)\!)$ whenever the Lie algebra 
$\mathfrak{g}$ is clear from the context.
Invoking a result of Korogodski and Soibelman (\cite{KorSoi-1998aa}), we will now describe  all simple unitarizable $\mathcal{O}_q(G)$-modules. 

\textbf{Elementary simple unitarizable $\mathcal{O}_q(G)$-modules:} For $1 \leq i \leq n$, let $d_i=\left\langle \alpha_i, \alpha_i \right\rangle /2$ and $q_i=q^{d_i}$.  
Define  $\phi_i :U_{q_i}(\mathfrak{sl}(2)) \longrightarrow  U_q(\mathfrak{g})$
be a  $*$-homomorphism given on the generators of $U_{q_i}(\mathfrak{sl}(2))$ by, 
\begin{displaymath}
 K \longmapsto K_i, \qquad E \longmapsto E_i, \qquad F \longmapsto F_i.
\end{displaymath}
By duality, it induces an  epimorphism 
\begin{displaymath}
\phi_i^* : \mathcal{O}_q(G) \longrightarrow \mathcal{O}_{q_i}(SU(2)).
\end{displaymath}
We  will use this map to get all elementary simple unitarizable modules of $\mathcal{O}_q(G)$.
Denote by $\Psi$ the following 
action of $\mathcal{O}_q(SU(2))$ on $c_{00}(\bbn)$ (see Proposition $4.1.1$, \cite{KorSoi-1998aa});
\begin{IEEEeqnarray}{rCl} \label{su(2)}
\Psi(u_l^k)e_p=\begin{cases}
              \sqrt{1-q^{2p}}e_{p-1} & \mbox{ if } k=l=1,\cr
              \sqrt{1-q^{2p+2}}e_{p+1} & \mbox{ if } k=l=2,\cr
							-q^{p+1}e_p & \mbox{ if } k=1,l=2,\cr
							q^pe_p & \mbox{ if } k=2,l=1,\cr
							\delta_{kl}e_p & \mbox{ otherwise }.\cr
              \end{cases}
\end{IEEEeqnarray}
For each $1 \leq i \leq n$, define an action  $\pi_{s_{i}}^n := \Psi \circ \phi_{i}^{*}$ of $\mathcal{O}_q(G)$. 
Each $\pi_{s_{i}}^n$ gives rise to 
an elementary simple $\mathcal{O}_q(G)$-module $V_{s_i}$. Also, for each $t \in \bbbt^n$, there are one dimensional  $\mathcal{O}_q(G)$-module $V_t$ with the action $\{\tau_t^n\}$.  
Given two actions  $\varphi$ and $\psi$ of $\mathcal{O}_q(G)$, define an action $\varphi * \psi := (\varphi \otimes \psi)\circ \Delta$. Similarly for any two $\mathcal{O}_q(G)$-module $V_{\varphi}$ and $V_{\phi}$, define 
$V_{\varphi} \otimes V_{\phi}$ as $\mathcal{O}_q(G)$-module with $\mathcal{O}_q(G)$ action coming from  $\varphi * \psi$.
For $w \in W_n$ such that  $s_{i_{1}}s_{i_{2}}...s_{i_{k}}$ is a reduced expression for $w$, define an action
$\pi_{w}^n$ as  $\pi_{s_{i_{1}}}^n*\pi_{s_{i_{2}}}^n*\cdots *\pi_{s_{i_{k}}}^n$  and  the corresponding $\mathcal{O}_q(G)$-module as $V_{w}$.  
Then $V_{w}$ is an simple unitarizable which is independent of the reduced expression. 
 Moreover, for $t \in \bbbt^{n}, w \in W_n$,  define an action $\pi_{t,w}^n$ by $\tau_{t}^n*\pi_{w}^n$ and denote the corresponding $\mathcal{O}_q(G)$-module by $V_{t,w}$. 
 We refer to (\cite{KorSoi-1998aa}, page~121) for the following theorem.
\bthm \label{allsimple}
The set $\left\{V_{t,w}; t \in \bbbt^{n}, w \in W_n\right\}$ is a complete set of mutually inequivalent simple  unitarizable left $\mathcal{O}_q(G)$-module.
\ethm

\newsubsection{Diagram representation}
As we are interested in the Gelfand-Kirillov dimension of certain modules it is essential that we understand the algebra action clearly. 
A diagram representation 
of various endomorphisms helps us in this respect. Corresponding to every $\pi_{s_i}^n$,  we have a diagram and thus a diagram for each $\pi_{w}^n$ 
is obtained as a concatenation of  diagrams of elementary modules. 
Using equation (\ref{su(2)}), one can write down explicitly  the endomorphisms  $\{\pi_{s_{i}}^n(u_l^k): 1 \leq k,l \leq N_n, 1 \leq i \leq n\}$  for  $G$ of type $A_n$, $C_n$ and $D_n$. 
We  order the indices of the simple roots in $\Pi$  in the standard way (see Appendix $C$ in \cite{Kir-2008aa}).
In this set up,
 these endomorphisms and their diagrams are given explicitly  in  \cite{ChaPal-2003ac} and   \cite{Sau-2015aa} for $G$ of type $A_n$ and $C_n$ respectively. 
For type $D_n$, it is given below. \\
For $i=1,2,\cdots,n-1$,
\[
\pi_{s_{i}}^n(u_l^k)e_p=\begin{cases}
              \sqrt{1-q^{2n}}e_{p-1} & \mbox{ if } (k,l)=(i,i) \mbox{ or } (2n-i,2n-i),\cr
            \sqrt{1-q^{2n+2}}e_{p+1} & \mbox{ if } (k,l)=(i+1,i+1) \mbox{ or } (2n-i+1,2n-i+1),\cr
							-q^{p+1}e_p & \mbox{ if } (k,l)=(i,i+1),\cr
							q^pe_p & \mbox{ if } (k,l)=(i+1,i),\cr
							q^{p+1}e_p & \mbox{ if } (k,l)=(2n-i,2n-i+1),\cr
							-q^p e_p& \mbox{ if } (k,l)=(2n-i+1,2n-i),\cr
							\delta_{kl}e_p & \mbox{ otherwise }. \cr
							\end{cases}
\]
For $i=n$,

\[
\pi_{s_n}^n(u_l^k)e_p=\begin{cases}
               \sqrt{1-q^{2n}}e_{p-1} & \mbox{ if } (k,l)=(n,n) \mbox{ or } (n-1,n-1),\cr
              \sqrt{1-q^{2n+2}}e_{p+1} & \mbox{ if } (k,l)=(n+1,n+1) \mbox{ or } (n+2,n+2),\cr
							-q^{p+1}e_p & \mbox{ if } (k,l)=(n-1,n+1),\cr
							q^pe_p & \mbox{ if } (k,l)=(n+1,n-1),\cr
							q^{p+1}e_p & \mbox{ if } (k,l)=(n,n+2)\cr
							-q^pe_p & \mbox{ if } (k,l)=(n+2,n),\cr
							\delta_{kl}e_p & \mbox{ otherwise }. \cr
              \end{cases}
\]
For $t=(t_{1},t_{2},\cdots ,t_{n}) \in \bbbt^{n}$,  the one dimensional irreducible representation  $\tau_t$ maps $u_j^i$ to $\overline{t_{i}}\delta_{ij}$ if $i \leq n$  and 
to $t_{2n+1-i}\delta_{ij}$ if $i > n$.

Now we will associate a diagram with each of the  above action of $\mathcal{O}_q(G)$. Let us describe the case of $\mathcal{O}_q(G)$ first.
For convenience, we use labeled lines and arrows to 
represent endomorphisms of a vector space  as given in the following table.
 \begin{center}
         \begin{tabular}{|c|c|c|c|}
   \hline
   Arrow type & Endomorphism &  Arrow type  & Endomorphism \\
   \hline   
\begin{tikzpicture}[scale=.7]
\draw [-] (0,0) -- (1,0);
\end{tikzpicture}
  &
  $I$ & 
  \begin{tikzpicture}[scale=.7]
\draw [-] (0,0) -- (1,0);
\node at (.7,.3) {$t$};
\end{tikzpicture}
& $M_t$  \\ 
     \hline   
 \begin{tikzpicture}[scale=.7]
 \draw [-] (0,0) -- (1,0);
 \node at (.7,.3) {$+$};
 \end{tikzpicture}
   &
   $S^*\sqrt{I-q^{2N+2}}$ & 
\begin{tikzpicture}[scale=.7]
\draw [-] (0,0) -- (1,0);
\node at (.5,-.2) {$-$};
\end{tikzpicture}
&
  $\sqrt{I-q^{2N+2}}S$  \\
     \hline
    \begin{tikzpicture}[scale=.7]
\draw [-] (0,.7) -- (.7,0);
\end{tikzpicture}
       &
$q^{N}$ & 
\begin{tikzpicture}[scale=.7]
\draw [-] (.1,.1) -- (.8,.8);
\end{tikzpicture}
  &
$-q^{N+1}$
\\
     \hline
     \begin{tikzpicture}[scale=.7]
\draw [->, very thick] (0,.7) -- (.7,0);
\end{tikzpicture}
       &
$-q^{N}$ & 
\begin{tikzpicture}[scale=.7]
\draw [->, very thick] (.1,.1) -- (.8,.8);
\end{tikzpicture}
  &
$q^{N+1}$
\\
     \hline
\end{tabular} \\ 
\end{center}

Note that for $t \in \bbbt$,  $M_t$ represents the endomorphism   on $\bbc$ sending $v$ to $tv$.   
For endomorphism $I$, the underlying vector space can be either $\bbc$ or $c_{00}(\bbn)$. For other endomorphisms, the vector spaces
on which they act are $c_{00}(\bbn)$. 
Let us describe how to use a diagram to represent the actions $\pi_{s_i}^n$ and $\tau_t^n$ where $1\leq i \leq n$ and $t \in \bbbt^{n}$.\\ 
 \begin{tabular}{p{150pt}p{150pt}p{150pt}}
\begin{tikzpicture}[scale=1.0, shift={(0,2)}]
\draw [-] (0,4.5) -- (1,4.5);
\draw [dashed] (0,4) -- (1,4);
\draw [-] (0,3.5) -- (1,3.5);
\draw [->, very thick] (0,3) -- (1,3.5);
\draw [->, very thick] (0,3.5) -- (1,3);
\draw [-] (0,3) -- (1,3);
\draw [dashed] (0,2.5) -- (1,2.5);
\draw [-] (0,2) -- (1,2);
\draw [-] (0,1.5) -- (1,1.5);
\draw [-] (0,2) -- (1,1.5);
\draw [-] (0,1.5) -- (1,2);
\draw [dashed] (0,1) -- (1,1);
\draw [-] (0,.5) -- (1,.5);
\node at (.5,3.7) {+};
\node at (.5,2.8) {-};
\node at (.5,2.2) {+};
\node at (.5,1.3) {-};
\node at (.5,5){$c_{00}(\bbn)$};
\node at (-.3,4.55){$2n$};
\node at (1.3,4.55){$2n$};
\node at (-1,3.55){$2n-i+1$};
\node at (2,3.55){$2n-i+1$};
\node at (-.9,3.05){$2n-i$};
\node at (1.9,3.05){$2n-i$};
\node at (-.5,2.05){$i+1$};
\node at (1.5,2.05){$i+1$};
\node at (-.3,1.55){$i$};
\node at (1.3,1.55){$i$};
\node at (-.3,.55){$1$};
\node at (1.3,.55){$1$};
 \node at (.5,0){\text Diagram 1: $\pi_{s_i}^n, i\neq n$};
\end{tikzpicture}
&
\begin{tikzpicture}[scale=1.0, shift={(0,2)}]
\draw [-] (0,4.5) -- (1,4.5);
\draw [dashed] (0,4) -- (1,4);
\draw [-] (0,3.5) -- (1,3.5);
\draw [-] (0,3) -- (1,2);
\draw [->, very thick] (0,3.5) -- (1,2.5);
\draw [->, very thick] (0,2.5) -- (1,3.5);
\draw [-] (0,2) -- (1,3);
\draw [-] (0,3) -- (1,3);
\draw [-] (0,2.5) -- (1,2.5);
\draw [-] (0,2) -- (1,2);
\draw [dashed] (0,1.5) -- (1,1.5);
\draw [-] (0,1) -- (1,1);
\node at (.5,3.7) {+};
\node at (.5,3.2) {+};
\node at (.5,2.3) {-};
\node at (.5,1.8) {-};
\node at (.5,5){$c_{00}(\bbn)$};
\node at (-.3,4.55){$2n$};
\node at (1.3,4.55){$2n$};
\node at (-.6,3.55){$n+2$};
\node at (1.6,3.55){$n+2$};
\node at (-.6,3.05){$n+1$};
\node at (1.6,3.05){$n+1$};
\node at (-.5,2.55){$n$};
\node at (1.5,2.55){$n$};
\node at (-.6,2.05){$n-1$};
\node at (1.6,2.05){$n-1$};
\node at (-.3,1.05){$1$};
\node at (1.3,1.05){$1$};
 \node at (.5,0){\text Diagram 2: $\pi_{s_n}^n$};
\end{tikzpicture}
&
\begin{tikzpicture}[scale=1.0, shift={(0,2)}]
\draw [-] (0,4.5) -- (1,4.5);
\draw [dashed] (0,4) -- (1,4);
\draw [-] (0,3.5) -- (1,3.5);
\draw [-] (0,3) -- (1,3);
\draw [-] (0,2.5) -- (1,2.5);
\draw [-] (0,2) -- (1,2);
\draw [dashed] (0,1.5) -- (1,1.5);
\draw [-] (0,1) -- (1,1);
\node at (.5,5.2){$\bbc$};
\node at (-.3,4.55){$2n$};
\node at (1.3,4.55){$2n$};
\node at (-.6,3.55){$n+2$};
\node at (1.6,3.55){$n+2$};
\node at (-.6,3.05){$n+1$};
\node at (1.6,3.05){$n+1$};
\node at (-.5,2.55){$n$};
\node at (1.5,2.55){$n$};
\node at (-.6,2.05){$n-1$};
\node at (1.6,2.05){$n-1$};
\node at (-.3,1.05){$1$};
\node at (1.3,1.05){$1$};
\node at (.5,4.7){$t_1$};
\node at (.5,3.7){$t_{n-1}$};
\node at (.5,3.2){$t_n$};
\node at (.5,2.7){$\overline{t_n}$};
\node at (.5,2.2){$\overline{t_{n-1}}$};
\node at (.5,1.2){$\overline{t_1}$};
 \node at (.5,0){\text Diagram 3: $\tau_t^n$};
\end{tikzpicture}

\end{tabular}\\[1ex]		
In these  diagrams, each path from a node $k$ on the
left to a node $l$ on the right stands for an
endomorphism   given as in the table acting on the vector space given at the top of the diagrams. Now 
$\pi_{s_i}^n(u_{l}^{k})$, $\pi_{s_n}^n(u_{l}^{k})$  and $\tau_{t}^n(u_{l}^{k})$ are the endomorphisms represented by the paths from $k$ to $l$ in 
diagram 1, diagram 2 and diagram 3 respectively
and are zero if there is no such path.
Thus, for example, in diagram 1, $\pi_{s_i}^n(u_1^1)$
is $I$; $\pi_{s_i}^n(u_1^2)$ is zero and $\pi_{s_i}^n(u_{i+1}^{i})=-q^{N+1}$ if $i >1$.

In case of $\mathcal{O}_q(SU(n+1))$, the table of symbols representing the operators appearing in the diagram of
$\pi_{s_i}^n$ for $1\leq i \leq n$ are as follows (see \cite{ChaPal-2003ac}).
\begin{center}
         \begin{tabular}{|c|c|c|c|}
   \hline
   Arrow type & Operator &  Arrow type  & Operator \\
   \hline   
\begin{tikzpicture}[scale=.7]
\draw [-] (0,0) -- (1,0);
\end{tikzpicture}
  &
  $I$ & 
  \begin{tikzpicture}[scale=.7]
\draw [-] (0,0) -- (1,0);
\node at (.7,.3) {$t$};
\end{tikzpicture}
& $M_t$  \\ 
     \hline   
 \begin{tikzpicture}[scale=.7]
 \draw [-] (0,0) -- (1,0);
 \node at (.7,.3) {$+$};
 \end{tikzpicture}
   &
   $S^*\sqrt{I-q^{2N+2}}$ & 
\begin{tikzpicture}[scale=.7]
\draw [-] (0,0) -- (1,0);
\node at (.5,-.2) {$-$};
\end{tikzpicture}
&
  $\sqrt{I-q^{2N+2}}S$  \\
     \hline
    \begin{tikzpicture}[scale=.7]
\draw [-] (0,.7) -- (.7,0);
\end{tikzpicture}
       &
$q^{N}$ & 
\begin{tikzpicture}[scale=.7]
\draw [-] (.1,.1) -- (.8,.8);
\end{tikzpicture}
  &
$-q^{N+1}$
\\
     \hline
\end{tabular} \\ 
\end{center}

 \begin{tabular}{p{200pt}p{200pt}}
\begin{tikzpicture}[scale=1.2, shift={(0,2)}]
\draw [-] (0,3) -- (1,3);
\draw [dashed] (0,2.5) -- (1,2.5);
\draw [-] (0,2) -- (1,2);
\draw [-] (0,2) -- (1,1);
\draw [-] (0,1) -- (1,2);
\draw [-] (0,1) -- (1,1);
\draw [dashed] (0,.5) -- (1,.5);
\draw [-] (0,0) -- (1,0);
\node at (.5,2.2) {+};
\node at (1,1.7) {};
\node at (1,1.3) {};
\node at (.5,.85) {${}-{}$};
\node at (-.5,3.05){$n+1$};
\node at (1.5,3.05){$n+1$};
\node at (-.5,2.05){$i+1$};
\node at (1.5,2.05){$i+1$};
\node at (-.3,1.05){$i$};
\node at (1.3,1.05){$i$};
\node at (-.3,.05){$1$};
\node at (1.3,.05){$1$};
\node at (.5,4){$L_2(\bbn)$};
 \node at (.5,-1){\text Diagram 4: $\pi_{s_i}$};
\end{tikzpicture}
&
\begin{tikzpicture}[scale=1.2, shift={(0,2)}]
 \draw [-] (0,3) -- (.5,3) node [above]{${}t_1$} -- (1,3);
 \draw [-] (0,2.5) -- (.5,2.5) node [above]{${}t_2$} -- (1,2.5);
 \draw [-] (0,2) -- (.5,2) node [above]{${}t_3$} -- (1,2);
  \draw [dashed] (0,1.5) --  (1,1.5);
 \draw [dashed] (0,1) -- (1,1);
\draw [-] (0,.5) -- (.5,.5) node [above]{${}t_{n}$} -- (1,.5);
 \draw [-] (0,0) -- (.5,0) node [above]{${}\overline{t_1}\cdots \overline{t_{n}}$} -- (1,0);
 \node at (-.5,3){$n+1$};
 \node at (-.5,2.5){$n$};
 \node at (-.5,2){$n-1$};
  \node at (-.5,.5){$2$};
  \node at (-.5,0){$1$};
 \node at (1.5,3){$n+1$};
 \node at (1.5,2.5){$n$};
 \node at (1.5,2){$n-1$};
  \node at (1.5,.5){$2$};
  \node at (1.5,0){$1$};
 \node at (.5,4){$\bbc$};
  \node at (.5,-1){\text Diagram 5: $\tau_{t}^n$};
\end{tikzpicture}
\end{tabular}\\[1ex]
For $\mathcal{O}_q(SP(2n))$, these data are described below (see \cite{Sau-2015aa}).\\
\begin{center}
         \begin{tabular}{|c|c|c|c|c|c|}
   \hline
   Arrow type & Endomorphism &  Arrow type  & Endomorphism & Arrow type  & Endomorphism\\
   \hline   
\begin{tikzpicture}[scale=.7]
\draw [-] (0,0) -- (1,0);
\end{tikzpicture}
  &
  $I$ & 
  \begin{tikzpicture}[scale=.7]
\draw [-] (0,0) -- (1,0);
\node at (.7,.3) {$t$};
\end{tikzpicture}
& $M_t$ & 
\begin{tikzpicture}[scale=.7]
 \draw [->>] (0,0) -- (1,.7);
 \end{tikzpicture}&  $-q^{2N+2}$  \\ 
     \hline   
 \begin{tikzpicture}[scale=.7]
 \draw [->] (0,0) -- (1,0);
 \node at (.7,.3) {$+$};
 \end{tikzpicture}
   &
   $S^*\sqrt{I-q^{2N+2}}$ & 
\begin{tikzpicture}[scale=.7]
\draw [-] (0,0) -- (1,0);
\node at (.5,-.2) {$-$};
\end{tikzpicture}
&
  $\sqrt{I-q^{2N+2}}S$ & 
\begin{tikzpicture}[scale=.7]
 \draw [->>] (0,0) -- (1,-.7);
 \end{tikzpicture}&  $q^{2N}$  \\
     \hline
    \begin{tikzpicture}[scale=.7]
\draw [-] (0,.7) -- (.7,0);
\end{tikzpicture}
       &
$q^{N}$ & 
\begin{tikzpicture}[scale=.7]
\draw [-] (.1,.1) -- (.8,.8);
\end{tikzpicture}
  &
$-q^{N+1}$ &&
\\
     \hline
     \begin{tikzpicture}[scale=.7]
\draw [->, very thick] (0,.7) -- (.7,0);
\end{tikzpicture}
       &
$-q^{N}$ & 
\begin{tikzpicture}[scale=.7]
\draw [->, very thick] (.1,.1) -- (.8,.8);
\end{tikzpicture}
  &
$q^{N+1}$ &&
\\
     \hline
\end{tabular} \\ 
\end{center}
\begin{tabular}{p{150pt}p{150pt}p{150pt}}

\begin{tikzpicture}[scale=1.0, shift={(0,2)}]
\draw [-] (0,4.5) -- (1,4.5);
\draw [dashed] (0,4) -- (1,4);
\draw [-] (0,3.5) -- (1,3.5);
\draw [->, very thick] (0,3) -- (1,3.5);
\draw [->, very thick] (0,3.5) -- (1,3);
\draw [-] (0,3) -- (1,3);
\draw [dashed] (0,2.5) -- (1,2.5);
\draw [-] (0,2) -- (1,2);
\draw [-] (0,1.5) -- (1,1.5);
\draw [-] (0,2) -- (1,1.5);
\draw [-] (0,1.5) -- (1,2);
\draw [dashed] (0,1) -- (1,1);
\draw [-] (0,.5) -- (1,.5);
\node at (.5,3.7) {+};
\node at (.5,2.8) {-};
\node at (.5,2.2) {+};
\node at (.5,1.3) {-};
\node at (.5,5){$c_{00}(\bbn)$};
\node at (-.3,4.55){$2n$};
\node at (1.3,4.55){$2n$};
\node at (-1,3.55){$2n-i+1$};
\node at (2,3.55){$2n-i+1$};
\node at (-.9,3.05){$2n-i$};
\node at (1.9,3.05){$2n-i$};
\node at (-.5,2.05){$i+1$};
\node at (1.5,2.05){$i+1$};
\node at (-.3,1.55){$i$};
\node at (1.3,1.55){$i$};
\node at (-.3,.55){$1$};
\node at (1.3,.55){$1$};
 \node at (.5,0){\text Diagram 6: $\pi_{s_i}^n, i\neq n$};
\end{tikzpicture}
&
\begin{tikzpicture}[scale=1.0, shift={(0,2)}]
\draw [-] (0,4.5) -- (1,4.5);
\draw [dashed] (0,4) -- (1,4);
\draw [-] (0,3.5) -- (1,3.5);
\draw [->>] (0,2.5) -- (1,3);
\draw [->>] (0,3) -- (1,2.5);
\draw [-] (0,3) -- (1,3);
\draw [-] (0,2.5) -- (1,2.5);
\draw [-] (0,2) -- (1,2);
\draw [dashed] (0,1.5) -- (1,1.5);
\draw [-] (0,1) -- (1,1);
\node at (.5,3.2) {+};
\node at (.5,2.3) {-};
\node at (.5,5){$c_{00}(\bbn)$};
\node at (-.3,4.55){$2n$};
\node at (1.3,4.55){$2n$};
\node at (-.6,3.55){$n+2$};
\node at (1.6,3.55){$n+2$};
\node at (-.6,3.05){$n+1$};
\node at (1.6,3.05){$n+1$};
\node at (-.5,2.55){$n$};
\node at (1.5,2.55){$n$};
\node at (-.6,2.05){$n-1$};
\node at (1.6,2.05){$n-1$};
\node at (-.3,1.05){$1$};
\node at (1.3,1.05){$1$};
 \node at (.5,0){\text Diagram 7: $\pi_{s_n}^n$};
\end{tikzpicture}
&
\begin{tikzpicture}[scale=1.0, shift={(0,2)}]
\draw [-] (0,4.5) -- (1,4.5);
\draw [dashed] (0,4) -- (1,4);
\draw [-] (0,3.5) -- (1,3.5);
\draw [-] (0,3) -- (1,3);
\draw [-] (0,2.5) -- (1,2.5);
\draw [-] (0,2) -- (1,2);
\draw [dashed] (0,1.5) -- (1,1.5);
\draw [-] (0,1) -- (1,1);
\node at (.5,5.2){$\bbc$};
\node at (-.3,4.55){$2n$};
\node at (1.3,4.55){$2n$};
\node at (-.6,3.55){$n+2$};
\node at (1.6,3.55){$n+2$};
\node at (-.6,3.05){$n+1$};
\node at (1.6,3.05){$n+1$};
\node at (-.5,2.55){$n$};
\node at (1.5,2.55){$n$};
\node at (-.6,2.05){$n-1$};
\node at (1.6,2.05){$n-1$};
\node at (-.3,1.05){$1$};
\node at (1.3,1.05){$1$};
\node at (.5,4.7){$t_1$};
\node at (.5,3.7){$t_{n-1}$};
\node at (.5,3.2){$t_n$};
\node at (.5,2.7){$\overline{t_n}$};
\node at (.5,2.2){$\overline{t_{n-1}}$};
\node at (.5,1.2){$\overline{t_1}$};
 \node at (.5,0){\text Diagram 8: $\tau_t^n$};
\end{tikzpicture}

\end{tabular}\\[1ex]	
Next, let us explain how to represent $\pi \ast \rho$ by a diagram where $\pi$ and $\rho$ are two actions of   $\mathcal{O}_q(G)$ acting on the vector  spaces $V_1$ and $V_2$ respectively.
Simply keep the two diagrams representing $\pi$ and $\rho$
adjacent to each other. Identify, for each row, the node on the right side
of the diagram for $\pi$ with the corresponding node on the left in the diagram
for $\rho$. Now, $(\pi \ast \rho)(u_{l}^{k})$ would be
an endomorphism of the vector space $V_1\otimes V_2$ 
determined by all the paths from
the node $k$ on the left to the node $l$ on the right. It would be zero if
there is no such path and if there are more than one paths, then it would be the sum of
the endomorphisms given by each such path.
In this way, we can draw diagrams for each action  of   $\mathcal{O}_q(G)$.

The following diagram  is for the action $\pi_{t,w}^4$ of   $\mathcal{O}_q(SO(8))$ where $w=s_1s_2s_3s_4s_2$ and $t=(t_1,t_2,1,1)$.

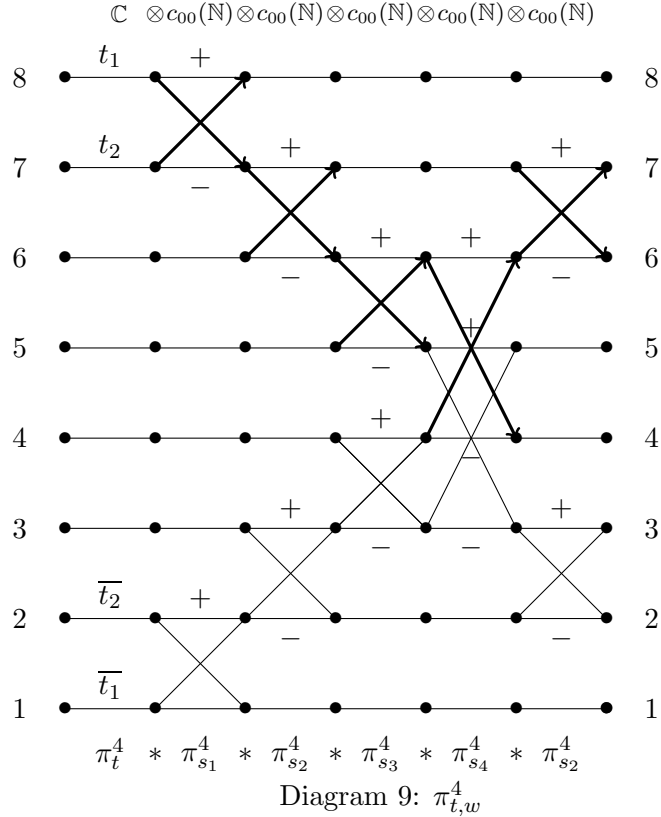
\begin{center}
 \begin{tikzpicture}[scale=1.2]
 \draw [-] (0,6) -- (.5,6) node [above]{${}t_1$} -- (1,6) -- (1.5,6) node [above]{${}+$} -- (2,6) -- (3,6) -- (4,6)  -- (5,6) -- (6,6);
 \draw [-] ((0,5) -- (.5,5) node [above]{${}t_2$} -- (1,5) -- (1.5,5) node [below]{${}-$} -- (2,5) -- (2.5,5) node [above]{${}+$} -- (3,5) -- (4,5)  -- (5,5) -- (5.5,5) node [above]{${}+$} -- (6,5);
 \draw [-] (0,4) -- (1,4) -- (2,4) -- (2.5,4) node [below]{${}-$} -- (3,4) -- (3.5,4)  node [above]{${}+$} -- (4,4)  -- (4.5,4)  node [above]{${}+$} -- (5,4) -- (5.5,4)  node [below]{${}-$} -- (6,4);
 \draw [-] (0,3) -- (1,3) -- (2,3) -- (3,3) -- (3.5,3) node [below]{${}-$} -- (4,3)  -- (4.5,3)  node [above]{${}+$} -- (5,3) -- (6,3);
 \draw [-] (0,2) -- (1,2) -- (2,2) -- (3,2) -- (3.5,2) node [above]{${}+$} -- (4,2) -- (4.5,2)  node [below]{${}-$} --   (5,2) -- (6,2);
 \draw [-] (0,1) -- (1,1)  -- (2,1) -- (2.5,1) node [above]{${}+$} -- (3,1) -- (3.5,1) node [below]{${}-$} --  (4,1) -- (4.5,1)  node [below]{${}-$} -- (5,1) -- (5.5,1)  node [above]{${}+$} -- (6,1);
 \draw [-] (0,0) -- (.5,0) node [above]{${}\overline{t_2}$} -- (1,0) -- (1.5,0) node [above]{${}+$} -- (2,0) -- (2.5,0) node [below]{${}-$} --  (3,0) -- (4,0) --   (5,0) -- (5.5,0)  node [below]{${}-$} -- (6,0);
 \draw [-] (0,-1) -- (.5,-1) node [above]{${}\overline{t_1}$} -- (1,-1)  -- (2,-1) --   (3,-1) --  (4,-1) --   (5,-1) -- (6,-1);
 \draw [->, very thick] (1,6) -- (2,5) node [below]{${}\;$};
 \draw [->, very thick] (2,5) -- (3,4) node [below]{${}\;$};
 \draw [->, very thick] (3,4) -- (4,3) node [below]{${}\;$};
 \draw [->, very thick] (4,4) -- (5,2); 
 \draw [-] (4,3) -- (5,1) node [below]{${}\;$};
 \draw [-] (3,2) -- (4,1) node [below]{${}\;$};
 \draw [-] (2,1) -- (3,0) node [below]{${}\;$};
 \draw [-] (1,0) -- (2,-1) node [below]{${}\;$};
 \draw [-] (5,1) -- (6,0) node [below]{${}\;$};
 \draw [->, very thick] (5,5) -- (6,4) node [below]{${}\;$};
 \draw [-] (3,2) -- (4,1) node [below]{${}\;$};
 \draw [->, very thick] (1,5) -- (2,6)  node [above]{${}\;$};
 \draw [->, very thick] (2,4) -- (3,5)  node [above]{${}\;$};
 \draw [->, very thick] (3,3) -- (4,4)  node [above]{${}\;$};
 \draw [->, very thick] (4,2) -- (5,4)  node [above]{${}\;$};
 \draw [-] (4,1) -- (5,3)  node [above]{${}\;$};
 \draw [-] (3,1) -- (4,2) node [above]{${}\;$};
 \draw [-] (2,0) -- (3,1)  node [above]{${}\;$};
 \draw [-] (1,-1) -- (2,0) node [above]{${}\;$};
 \draw [-] (5,0) -- (6,1)  node [above]{${}\;$};
 \draw [->, very thick] (5,4) -- (6,5)  node [above]{${}\;$};
 \draw [-] (3,1) -- (4,2)  node [above]{${}\;$};
 \node at (0,6){$\bullet$};
 \node at (1,6){$\bullet$};
 \node at (2,6){$\bullet$};
 \node at (3,6){$\bullet$};
 \node at (4,6){$\bullet$};
 \node at (5,6){$\bullet$};
 \node at (6,6){$\bullet$};
 \node at (0,5){$\bullet$};
 \node at (1,5){$\bullet$};
 \node at (2,5){$\bullet$};
 \node at (3,5){$\bullet$};
 \node at (4,5){$\bullet$};
 \node at (5,5){$\bullet$};
 \node at (6,5){$\bullet$};
 \node at (0,4){$\bullet$};
 \node at (1,4){$\bullet$};
 \node at (2,4){$\bullet$};
 \node at (3,4){$\bullet$};
 \node at (4,4){$\bullet$};
 \node at (5,4){$\bullet$};
 \node at (6,4){$\bullet$};
 \node at (0,3){$\bullet$};
 \node at (1,3){$\bullet$};
 \node at (2,3){$\bullet$};
 \node at (3,3){$\bullet$};
 \node at (4,3){$\bullet$};
 \node at (5,3){$\bullet$};
 \node at (6,3){$\bullet$};
 \node at (0,2){$\bullet$};
 \node at (1,2){$\bullet$};
 \node at (2,2){$\bullet$};
 \node at (3,2){$\bullet$};
 \node at (4,2){$\bullet$};
 \node at (5,2){$\bullet$};
 \node at (6,2){$\bullet$};
 \node at (0,1){$\bullet$};
 \node at (1,1){$\bullet$};
 \node at (2,1){$\bullet$};
 \node at (3,1){$\bullet$};
 \node at (4,1){$\bullet$};
 \node at (5,1){$\bullet$};
 \node at (6,1){$\bullet$};
 \node at (0,0){$\bullet$};
 \node at (1,0){$\bullet$};
 \node at (2,0){$\bullet$};
 \node at (3,0){$\bullet$};
 \node at (4,0){$\bullet$};
 \node at (5,0){$\bullet$};
  \node at (6,0){$\bullet$};
 \node at (6,-1){$\bullet$};
 \node at (0,-1){$\bullet$};
 \node at (1,-1){$\bullet$};
 \node at (2,-1){$\bullet$};
 \node at (3,-1){$\bullet$};
 \node at (4,-1){$\bullet$};
 \node at (5,-1){$\bullet$};
 \node at (6,-1){$\bullet$};
 \node at (.6,6.7){\footnotesize{$\bbc$}};
 \node at (1.5,6.7){\footnotesize{$c_{00}(\bbn)$}};
 \node at (2.5,6.7){\footnotesize{$c_{00}(\bbn)$}};
 \node at (3.5,6.7){\footnotesize{$c_{00}(\bbn)$}};
 \node at (4.5,6.7){\footnotesize{$c_{00}(\bbn)$}};
 \node at (5.5,6.7){\footnotesize{$c_{00}(\bbn)$}};
 \node at (1,6.7){\footnotesize{$\otimes$}};
  \node at (2,6.7){\footnotesize{$\otimes$}};
  \node at (3,6.7){\footnotesize{$\otimes$}};
 \node at (4,6.7){\footnotesize{$\otimes$}};
  \node at (5,6.7){\footnotesize{$\otimes$}};
   \node at (-.5,6){$8$};
   \node at (-.5,5){$7$};
   \node at (-.5,4){$6$};
   \node at (-.5,3){$5$};
   \node at (-.5,2){$4$};
   \node at (-.5,1){$3$};
     \node at (-.5,0){$2$};
       \node at (-.5,-1){$1$};
   \node at (6.5,6){$8$};
   \node at (6.5,5){$7$};
   \node at (6.5,4){$6$};
   \node at (6.5,3){$5$};
   \node at (6.5,2){$4$};
   \node at (6.5,1){$3$};
    \node at (6.5,0){$2$};
       \node at (6.5,-1){$1$};
        \node at (.5,-1.5){$\pi_{t}^4$};
          \node at (1.5,-1.5){$\pi_{s_1}^4$};
            \node at (2.5,-1.5){$\pi_{s_2}^4$};
              \node at (3.5,-1.5){$\pi_{s_3}^4$};
                \node at (4.5,-1.5){$\pi_{s_4}^4$};
                  \node at (5.5,-1.5){$\pi_{s_2}^4$};
                   \node at (1,-1.5){$*$};
            \node at (2,-1.5){$*$};
              \node at (3,-1.5){$*$};
                \node at (4,-1.5){$*$};
                  \node at (5,-1.5){$*$};
 \node at (3.5,-2){\text Diagram 9: $\pi_{t,w}^4$};
 \end{tikzpicture}
 \end{center}
Here there are two paths from the vertex $1$  in the left to the vertex $2$ in the right
 and hence the endomorphism 
$\pi_{t,w}^4(u_1^2)$   will be sum of the  endomorphisms represented by each path. Therefore we have  
\begin{IEEEeqnarray}{rCl}
 \pi_{t,w}^4(u_1^2)&=&\overline{t_1}q^{N+1}\otimes q^{N+1}\otimes \sqrt{1-q^{2N+2}}S^*\otimes \sqrt{1-q^{2N+2}}S^*\otimes q^N \nonumber  \\
 &&-\overline{t_1}q^{N+1}\otimes  \sqrt{1-q^{2N+2}}S^*\otimes I \otimes I \otimes \sqrt{1-q^{2N+2}}S^*. \nonumber 
\end{IEEEeqnarray}

Let us record few properties of these endomorphisms which are immidiate from the diagram.
\begin{enumerate}
\item 
The endomorphism $\pi_w^4(u_6^8)$ when applied to vector of the form $e_{i_1} \otimes e_{i_2}\otimes e_{i_3} \otimes e_{i_4} \otimes e_0$ gives $Ce_{i_1} \otimes e_{{i_2}+1}\otimes e_{i_3} \otimes e_{i_4} \otimes e_0$. 
The point is; if  a standard basis element  has $e_0$ in last co-ordinate then there is only one path from the $8$ to $6$ that survives as $Se_0=0$ and it increases second co-ordinate of the standard basis element by a unit. 
 \item 
 For each vertex $i \in \{2,3,\cdots 7\}$ in the left, there is an unique vertex, say $r_w(i)$ in the right such that $\pi_w^4(u_{r(i)}^i)(e_0\otimes e_0\otimes \cdots \otimes e_0)=Ce_0\otimes e_0\otimes \cdots \otimes e_0$. 
More explicitly, the map $r_w$ sends $2 \mapsto 1, 3 \mapsto 3,  4 \mapsto 5,  5 \mapsto 4,  6 \mapsto 6$ and  $7 \mapsto 8$.
 \item 
 Also, there may be many paths from $i$ to $r_w(i)$  
 but  when applied to $e_0\otimes e_0\otimes \cdots \otimes e_0$, only one path give a nonzero vector i.e  $Ce_0\otimes e_0\otimes \cdots \otimes e_0$.
 \item 
For any $j \in \{2,3,\cdots 7\}$, $j \neq i$,  we have $\pi_w^4(u_{r_w(i)}^j)(e_0\otimes e_0\otimes \cdots \otimes e_0)=0$.
\end{enumerate} 
As we can see that there is nothing special about the above diagram and  these observations hold in general. 
We will make these things more precise in the next section but with this example in mind one can get a better understanding of the results. 

\newsubsection{Weyl group}
We will now give a brief introduction of  the Weyl group $W_n$ of type $A_n$, $C_n$ and $D_n$. For details, we refer the reader to \cite{FulHar-1991aa} or \cite{Hum-1990aa}. 
Let $E_{i,j}$ denote the $n \times n$ matrix having only one non-zero entry $1$ at $ij^{th}$ position. 
The group  $W_n$ is isomorphic to a subgroup of $GL(n,\bbr)$ generated by $s_{1},s_{2},...s_{n}$ where \\
\textbf{ Case 1}: $\mathfrak{sl}(n+1)$
\begin{IEEEeqnarray}{rCll}
s_{i} &=& I-E_{i,i}-E_{i+1,i+1}+E_{i,i+1}+E_{i+1,i}, & \qquad \mbox{for}\quad i=1,2,...n.\nonumber
\end{IEEEeqnarray}
Hence $W_n$ is isomorphic to permutation group of $n$ letters. \\
\textbf{ Case 2}: $\mathfrak{sp}(2n)$
\begin{IEEEeqnarray}{rCll}
s_{i} &=& I-E_{i,i}-E_{i+1,i+1}+E_{i,i+1}+E_{i+1,i}, & \qquad \mbox{for}\quad i=1,2,...n-1,\nonumber \\
s_{n} &=& I-2E_{n,n}, & \qquad \mbox{for} \quad i=n.\nonumber
\end{IEEEeqnarray}
Hence $W_n$ is the set of $n \times n$ matrices having one non-zero entry in each row and each column which is either $1$ or $-1$.\\
\textbf{ Case 3}: $\mathfrak{so}(2n)$
\begin{IEEEeqnarray}{rCll}
s_{i} &=& I-E_{i,i}-E_{i+1,i+1}+E_{i,i+1}+E_{i+1,i}, & \qquad \mbox{for}\quad i=1,2,...n-1,\nonumber \\
s_{n} &=& I-2E_{n,n}-2E_{n-1,n-1}, & \qquad \mbox{for} \quad i=n.\nonumber
\end{IEEEeqnarray}
Hence $W_n$ is the set of $n \times n$ matrices having one non-zero entry in each row and each column which is either $1$ or $-1$ and number of $-1$'s are even.

Any element of $W_n$ can be written in the form: 
$\psi_{1,k_{1}}^{(\epsilon_{1})}(w)\psi_{2,k_{2}}^{(\epsilon_{2})}(w)\cdots \psi_{n,k_{n}}^{(\epsilon_{n})}(w)$ for some choices 
of $\epsilon_1, \epsilon_2, \cdots, \epsilon_n$ and $k_1,k_2,\cdots k_n$ where
$\epsilon_{r} \in \left\{0,1,2\right\}$ and $n-r+1\leq k_{r}\leq n$ with the convention that, \\
\textbf{Case 1}: $\mathfrak{sl}(n+1)$
\[
 \psi_{r,k_{r}}^{\epsilon}(w)= \begin{cases}
               s_{r}s_{r-1}\cdots s_{n-k_r+1} & \mbox{ if } \epsilon=1,2\cr
							 \mbox{ empty string } & \mbox{ if } \epsilon=0.\cr
               \end{cases}
\]
\textbf{Case 2}: $\mathfrak{sp}(2n)$
\[
\psi_{r,k_{r}}^{\epsilon}(w)=\begin{cases}
               s_{n-r+1}s_{n-r+2}\cdots s_{k_r} & \mbox{ if } \epsilon=1,\cr
               s_{n-r+1}s_{n-r+2}\cdots...s_{n-1}s_{n}s_{n-1}\cdots s_{k_{r}} & \mbox{ if } \epsilon=2,\cr
							 \mbox{ empty string } & \mbox{ if } \epsilon=0.\cr
               \end{cases}
\]
\textbf{Case 3}: $\mathfrak{so}(2n)$
\[
\psi_{r,k_{r}}^{\epsilon}(w)=\begin{cases}
               s_{n-r+1}s_{n-r+2}\cdots s_{k_r} & \mbox{ if } \epsilon=1,\cr
               s_{n-r+1}s_{n-r+2}\cdots...s_{n-1}s_{n}s_{n-2}s_{n-3}\cdots s_{k_{r}} & \mbox{ if } \epsilon=2,\cr
							 \mbox{ empty string } & \mbox{ if } \epsilon=0.\cr
               \end{cases}
\]
Also, the above expression is a reduced expression. We will call the word $\psi_{r,k_{r}}^{(\epsilon_{r})}(w)$ 
the $r^{th}$ part $w_r$ of $w$. 
Accordingly, the simple module $V_w$ decomposes as $V_w=V_{w_1}\otimes V_{w_2}\otimes \cdots \otimes V_{w_n}$.
We call $V_{w_i}$ the $i^{th}$ part of $V_w$.


\newsection{Main result}
In this section, our main aim is to compute the Gelfand-Kirillov dimension of $V_{t,w}$ for $t \in \bbbt^n$ and $w \in W_n$. 
We first recall from \cite{KraLen-2000aa}
the definition of Gelfand Kirillov dimension of a module.
\bdfn (\cite{KraLen-2000aa})
Let $A$ be a  unital algebra and $M$ be a  left $A$-module. The Gelfand-Kirillov dimension of $M$ is given by 
\[
 \mbox{GKdim}(M)=sup_{V,F}\varlimsup\frac{\ln \dim(V^kF)}{\ln k}  
\]
where the supremum is taken over all finite dimensional subspace  $V$ of $A$ containing $1$  and all finite dimensional subspaces $F$ of $M$.  If 
 $A$ is a finitely generated unital algebra and $M$ be a finitely generated left $A$-module then
\[
 \mbox{GKdim}(M)=sup_{\xi,F}\varlimsup\frac{\ln \dim(\xi^kF)}{\ln k}  
\]
where the supremum is taken over all finite sets  $\xi$  containing $1$ that generates $A$ and all finite dimensional subspaces $F$ of $M$ that generates $M$ as an left $A$-module.
\edfn
\brmrk  
 For a finitely generated algebra $A$ and a finitely generated left $A$-module $M$, the quantity ``$\varlimsup\frac{\ln \dim(\xi^kF)}{\ln k}$'' does not depend on particular choices of $\xi$ and $F$. Therefore 
  one can choose a fixed (but finite) set of generators of $A$ and  
a finite dimensional subspaces $F$ of $M$ that generates $M$ as a left $A$-module. Moreover, if $M$ is a simple left $A$-module then one can take one dimensional subspace spanned by any vector $v \in M$ as 
a candidate for $F$ (see \cite{KraLen-2000aa} for details). 
\ermrk
Let $\xi_n=\{u_j^i:1\leq i,j \leq N_n\}$ where $N_n$ is the dimension of the representation $(\!(u_j^i)\!)$. Define
\begin{IEEEeqnarray}{rCl}
 \xi_{G_q}=\begin{cases}
            \xi_n \cup \{1\} & \mbox{ for } \mathcal{O}_q(G)=\mathcal{O}_q(SU{(n+1)}) \mbox{ or } \mathcal{O}_q(SP(2n)), \cr
            \xi_n \cup \{z_j^i:1\leq i,j \leq 2^n\} \cup \{1\} & \mbox{ for } \mathcal{O}_q(G)=\mathcal{O}_q(\mbox{Spin}(2n)). \cr
           \end{cases} \nonumber
\end{IEEEeqnarray} 
Then $\xi_{G_q}$ is a generating set of $\mathcal{O}_q(G)$ containing $1$. Throughout the article, we will work with this generating set. 
Let $w \in W_n$ and $t \in \bbbt^n$. 
The following lemma 
says that GKdim$(V_{w,t})$ is less than or equal to the length of the Weyl word.

\blmma \label{upper}
 Assume that $\mathcal{O}_q(G)$ is  one of the Hopf $*$-algebras $\mathcal{O}_q(SU(n+1))$, $\mathcal{O}_q(SP(2n))$ or $\mathcal{O}_q(\mbox{Spin}(2n))$.
For $w\in W_n$ and $t \in \bbbt^n$, let $V_{t,w}$ be the associated simple unitarizable left $\mathcal{O}_q(G)$-module. Then one has 
GKdim$(V_{t,w})\leq \ell(w)$.
\elmma
\prf 
 The algebra $\mathcal{O}_q(SU(2))$ has the following  vector space basis  (see Proposition 4,  page 100, \cite{KliSch-1997aa}).
 \[
  \{(u_{1,\mathfrak{sl}_2}^1)^b(u_{1,\mathfrak{sl}_2}^2)^c(u_{2,\mathfrak{sl}_2}^1)^d, 
 (u_{1,\mathfrak{sl}_2}^2)^v(u_{2,\mathfrak{sl}_1}^1)^w (u_{2,\mathfrak{sl}_2}^2)^y : b,c,d,v,w,y \in  \bbn\}.
 \]
Define the 
exponent of the element $(u_{1,\mathfrak{sl}_2}^1)^b(u_{1,\mathfrak{sl}_2}^2)^c(u_{2,\mathfrak{sl}_2}^1)^d$ to be $b$ and  the exponent of the element
$(u_{1,\mathfrak{sl}_2}^2)^v(u_{2,\mathfrak{sl}_1}^1)^w (u_{2,\mathfrak{sl}_2}^2)^y$ to be $y$.   Given a linear 
combination of these elements, define its exponent to be maximum of the exponent of the terms appearing in the linear combination.
Let
\[
 a=\max_{u \in \xi_{G_q}, i\in \{1,2,\cdots,n\}}\big\{\mbox{exponent of } \phi_i^*(u): u \in \xi_{G_q}\big\}.
\]
Using equation (\ref{su(2)}), we get 
\[
       \mbox{span }\pi_{t,w}^n(\xi_{G_q}^r)(e_0\otimes \cdots \otimes e_0)\subset \mbox{span}\{e_{i_1} \otimes \cdots \otimes e_{i_{\ell(w)}}: 0\leq i_j \leq ra \mbox{ for all } 1\leq j \leq \ell(w)\}.
      \]
Therefore one has
\[
 \mbox{GKdim}(V_{t,w})=\varlimsup\frac{\ln \dim(\xi_{G_q}^r(e_0\otimes \cdots \otimes e_0))}{\ln r} \leq \varlimsup\frac{\ln (ra)^{\ell(w)}}{\ln r}=\ell(w).
\]
This proves the claim.
\qed 

In the following lemma, we give some polynomials  which when applied in a certain order to a fixed vector of the form  $v \otimes \underbrace{e_0\otimes e_0\otimes \cdots \otimes e_0}_{n^{th}-\mbox{part}}$ 
give all matrix units of the $n^{th}$ part. More precisely,

  
\blmma \label{sample1}
 Assume that $\mathcal{O}_q(G)$ is  one of the Hopf $*$-algebras $\mathcal{O}_q(SU(n+1))$, $\mathcal{O}_q(SP(2n))$ or $\mathcal{O}_q(\mbox{Spin}(2n))$ 
 and  $w \in W_n$. 
Then there exist polynomials $p_{1}^{(w,n)}, p_{2}^{(w,n)},\cdots p_{\ell(w_n)}^{(w,n)}$ with noncommuting variables $\pi_w^n(u_j^{N_n})$'s 
and a permutation $\sigma$ of $\{1,2,\cdots ,\ell(w_n)\}$ such that  for all $r_1^n,r_2^n,\cdots,r_{\ell(w_n)}^n \in \bbn$, one has                                                                                                         
\begin{IEEEeqnarray}{lCl} 
 (p_{1}^{(w,n)})^{r_{\sigma(1)}^n}(p_{2}^{(w,n)})^{r_{\sigma(2)}^n}\cdots (p_{\ell(w_n)}^{(w,n)})^{r_{\sigma(\ell(w_n))}^n} (v \otimes e_0\otimes \cdots \otimes e_0) 
 = C v \otimes  e_{r_{1}^n}\otimes e_{r_{2}^n} \otimes \cdots \otimes e_{r_{\ell(w_n)}^n} \nonumber 
\end{IEEEeqnarray}
where $v \in c_{00}(\bbn)^{\otimes \sum_{i=1}^{n-1} \ell(w_i)}$ and  $C$ is a nonzero real number.
\elmma
\prf Let us consider each case separately. \\ 
\textbf{Case 1:}  $\mathcal{O}_q(G) =\mathcal{O}_q(SU(n+1))$.  In this case $N_n=n+1$.\\
  For $1\leq j \leq \ell(w_n)$, define $\sigma(j):=j$ and  $p_{j}^{(w,n)}:=\pi_w^n(u_{N_n-j+1}^{N_n})$.  
  Using diagram representations of $\mathcal{O}_q(SU(n+1))$, one has 
\begin{IEEEeqnarray}{rCl} \label{sun}
 p_{j}^{(w,n)}&=&1^{\otimes \sum_{i=1}^{n-1} \ell(w_i)} \otimes (q^{N})^{\otimes j-1} \otimes \underbrace{\sqrt{1-q^{2N}}S^*}_{j^{th} \mbox{ place}}\otimes 1^{\otimes \ell(w_n)-j}. \nonumber    
\end{IEEEeqnarray}  
Therefore
\begin{IEEEeqnarray}{lCl} 
 ( p_{j}^{(w,n)})^r&\sim & 1^{\otimes \sum_{i=1}^{n-1} \ell(w_i)} \otimes (q^{rN})^{\otimes j-1} \otimes \underbrace{(\sqrt{1-q^{2N}}S^*)^r}_{j^{th} \mbox{ place}}\otimes 1^{\otimes \ell(w_n)-j},  \nonumber\\
 &=& 1^{\otimes \sum_{i=1}^{n-1} \ell(w_i)} \otimes 1^{\otimes j-1} \otimes \underbrace{(\sqrt{1-q^{2N}}S^*)^r}_{j^{th} \mbox{ place}}\otimes 1^{\otimes \ell(w_n)-j}. \nonumber
\end{IEEEeqnarray}
on the whole  vector space. Using this, one can get the claim.\\
\textbf{Case 2 :} $\mathcal{O}_q(G) =\mathcal{O}_q(SP(2n))$. In this case $N_n=2n$.\ \\
If $\ell(w_n) \leq n$ then $w=w_1w_2\cdots w_{n-1}s_1s_{2}\cdots s_{\ell(w_n)}$. In this case, define  $\sigma(j):=j$ and  $p_{j}^{(w,n)}:=\pi_w^n(u_{N_n-j+1}^{N_n})$. If $\ell(w_n)> n$, then 
$w=w_1w_2\cdots w_{n-1}s_1s_2\cdots s_{n-1}s_n s_{n-1}\cdots s_{N_n-\ell(w_n)}$. In this case, define
\begin{IEEEeqnarray}{lCl}
 p_{j}^{(w,n)}= \begin{cases}
                   \pi_w^n(u_{N_n-j+1}^{N_n}) & \mbox{ if }  1 \leq  j\leq  N_n-\ell(w_n) -1, \cr
                   \pi_w^n(u_{j+1}^{N_n})  & \mbox{ if } N_n-\ell(w_n) \leq j <n, \cr
                    [\pi_w^n(u_{n+1}^{N_n}),\pi_w^n(u_{n}^{N_n})]  & \mbox{ if }  j=n, \cr
                    \pi_w^n(u_{j}^{N_n}) & \mbox{ if } n+1 \leq j \leq \ell(w_n). \cr
                  \end{cases} \nonumber
\end{IEEEeqnarray} 
Let $\sigma(j)$ be $N_n-j$ if $N_n-\ell(w_n) \leq j \leq \ell(w_n)$ and $j$ otherwise. By a direct computation or using diagram representation, we have
\begin{IEEEeqnarray}{lCl} \label{sp2n}
  p_{j}^{(w,n)} \sim \begin{cases}
                         1^{\otimes \sum_{i=1}^{n-1} \ell(w_i)}\otimes 1^{\otimes \sigma(j)-1}  \otimes \underbrace{\sqrt{1-q^{2N}}S^*}_{\sigma(j)^{th} \mbox{ place}}\otimes 1^{\otimes \ell(w_n)- \sigma(j)}  & \mbox{ if } j \neq n, \cr
                       1^{\otimes \sum_{i=1}^{n-1} \ell(w_i)} \otimes 1^{\otimes \sigma(n)-1}  \otimes \underbrace{\sqrt{1-q^{4N}}S^*}_{\sigma(n)^{th} \mbox{ place}}\otimes 1^{\otimes \ell(w_n)- \sigma(n)}  & \mbox{ if } j=n, \cr
                       \end{cases}
\end{IEEEeqnarray}
on the subspace generated by standard  basis elements having $e_0$ at $\sigma(k)^{th}$ place for $k <j$. Since $q^N\sqrt{1-q^{2N}}S^*=\sqrt{1-q^{2N}}S^*q^{N+2}$ and 
$q^N\sqrt{1-q^{4N}}S^*=\sqrt{1-q^{4N}}S^*q^{N+4}$, we get
\begin{IEEEeqnarray}{lCl}
  (p_{j}^{(w,n)})^r \sim \begin{cases}
                         1^{\otimes \sum_{i=1}^{n-1} \ell(w_i)}\otimes 1^{\otimes \sigma(j)-1}  \otimes \underbrace{(\sqrt{1-q^{2N}}S^*)^r}_{\sigma(j)^{th} \mbox{ place}}\otimes 1^{\otimes \ell(w_n)- \sigma(j)}  & \mbox{ if } j \neq n, \cr
                         1^{\otimes \sum_{i=1}^{n-1} \ell(w_i)}\otimes 1^{\otimes \sigma(n)-1}  \otimes \underbrace{(\sqrt{1-q^{4N}}S^*)^r}_{\sigma(n)^{th} \mbox{ place}}\otimes 1^{\otimes \ell(w_n)- \sigma(n)}  & \mbox{ if } j=n, \cr
                       \end{cases} \nonumber
\end{IEEEeqnarray}
The rest is a straightforward checking. \\
\textbf{Case 3:} $\mathcal{O}_q(G) =\mathcal{O}_q(\mbox{Spin}(2n))$. In this case $N_n=2n$.\ \\
If $\ell(w_n) < n$ then $w=w_1w_2\cdots w_{n-1}s_1s_{2}\cdots s_{\ell(w_n)}$. In this case, define   $p_{j}^{(w,n)}:=\pi_w^n(u_{N_n-j+1}^{N_n})$ and $\sigma(j):=j$.
If $\ell(w_n) \geq n$, then 
$w=w_1w_2\cdots w_{n-1}s_1s_2\cdots s_{n-2}s_{n-1}s_n s_{n-2}\cdots s_{N_n-\ell(w_n)-1}$. In this case, define 
\begin{IEEEeqnarray}{rCl} \label{q}
 p_{j}^{(w,n)}= \begin{cases}
                   \pi_w^n(u_{N_n-j+1}^{N_n}) & \mbox{ if }    1 \leq  j\leq  N_n-\ell(w_n) -2, \cr
                   \pi_w^n(u_{j+1}^{N_n})  & \mbox{ otherwise}. \cr
                  \end{cases} 
\end{IEEEeqnarray} 
Let $\sigma(j)$ be $N_n-j-1$ if $ N_n-\ell(w_n)-1 \leq j \leq \ell(w_n)$ and $j$ otherwise. Using diagram representation, we have
\begin{IEEEeqnarray}{lCl} \label{so2n}
 p_{j}^{(w,n)} \sim 1^{\otimes \sum_{i=1}^{n-1} \ell(w_i)}\otimes 1^{\otimes \sigma(j)-1}  \otimes \underbrace{\sqrt{1-q^{2N}}S^*}_{\sigma(j)^{th} \mbox{ place}}\otimes 1^{\otimes \ell(w_n)- \sigma(j)} 
\end{IEEEeqnarray}
on the subspace generated by standard  basis elements having $e_0$ at $\sigma(k)^{th}$ place for $k <j$. The similar calculations 
as above will prove the result.
\qed
\brmrk
\begin{enumerate}
 \item As we saw in part $(1)$ of the observations made in subsection $2.1$ that the order in which we apply these endomorphisms are important just to make sure that the matrix unit $e_0$ 
 remains in certain components and that's why a permutation occurs in this lemma. 

 \item
Note that the polynomials $p_{j}^{(w,n)}$; $1 \leq j \leq \ell(w_n)$ defined here involve   elements of last row of the matrix $(\!(u_j^i)\!)$ only.
\end{enumerate}

\ermrk
To show our main claim, we need to extend these results  for all parts of the Weyl word. 
For type $C_n$ and $D_n$, let  $M_n^i=n-i+1$  and $N_n^i:=N_n-n+i$ and for type $A_n$, let
$M_n^i=1$ and $N_n^i:=i+1$.
\blmma \label{unique path} Assume that $\mathcal{O}_q(G)$ is  one of 
the Hopf $*$-algebras $\mathcal{O}_q(SU(n+1))$, $\mathcal{O}_q(SP(2n))$ or $\mathcal{O}_q(\mbox{Spin}(2n))$.
\begin{enumerate}
 \item Let $w \in W_n$ be of the form   $w=w_{i+1}$ and let  $V_w$ be the associated   $\mathcal{O}_q(G)$-module. 
 Then for each $M_n^i \leq k \leq N_n^i$, there exists unique $r_w(k) \in \{M_n^{i+1},\cdots ,N_n^{i+1}\}$ 
  such that 
\begin{IEEEeqnarray}{rCl} \label{m1}
 \pi_w^n(u_{r_w(k)}^k)(e_0\otimes e_0\otimes \cdots \otimes e_0)=Ce_0\otimes e_0\otimes \cdots \otimes e_0 
\end{IEEEeqnarray} 
 where $C$ is a nonzero real number. 
Moreover, 
\begin{IEEEeqnarray}{rCl}\label{m2}
\pi_w^n(u_{r_w(k)}^j)(e_0\otimes e_0\otimes \cdots \otimes e_0)=0 \mbox{ for } j \in \{M_n^i,\cdots ,N_n^i\}/\{k\}. 
\end{IEEEeqnarray}
\item
Let $w \in W_n$ be of the form   $w=w_{i+1}w_{i+2}\cdots w_l$, $l\leq n$ and let  $V_w$ be the associated   $\mathcal{O}_q(G)$-module. 
For each $M_n^i \leq k \leq N_n^i$, define $r_w(k):=r_{w_l} \circ r_{w_{l-1}} \circ \cdots r_{w_{i+2}}\circ r_{w_{i+1}}(k) \in \{M_n^{l},\cdots ,N_n^{l}\}$. Then $r_w$ 
satisfies equations (\ref{m1}) and (\ref{m2}).
\end{enumerate}

\elmma
We have seen these results  in a particular case given  as   part $(2),(3)$ and $(4)$ of the observations   in subsection $2.1$. This lemma basically  generalises these observations to an appropriate set up. \\
\prf \textbf{Proof of part (1):}
To show the claim, we will define the function $r_{w_{i+1}}$ explicitly in each case in the following way. \\ 
\textbf{Case 1:} $\mathcal{O}_q(G) =\mathcal{O}_q(SU(n+1))$. \\  In this case, $w_{i+1}=s_is_{i-1} \cdots s_{i-\ell(w_{i+1})-1}$ and hence for all $1\leq k \leq n$, $s_k$ occurs in $w_{i+1}$ either once or does not occur. 
\[
 r_{w_{i+1}}(j)=\begin{cases}
                j+1, & \mbox{ if } s_j \mbox{ occurs in } w_{i+1} \mbox{ once }, \cr
                j &  \mbox{ otherwise }. \cr                                                      
               \end{cases} 
               \]
\textbf{Case 2 :} $\mathcal{O}_q(G) =\mathcal{O}_q(SP(2n))$.\\ In this case, $w_{i+1}$ is either of the form $s_{i}s_{i+1}\cdots s_k$ where $k \leq n$ or   $s_{i}s_{i-1}\cdots s_{n-1}s_ns_{n-1}\cdots  s_k$. Therefore
for all $1\leq k < n$, $s_k$ occurs in  $w_{i+1}$ either once or  twice or does not occur. For $j \leq n$, define
\[
 r_{w_{i+1}}(j)=\begin{cases}
                j-1, & \mbox{ if } s_{j-1} \mbox{ occurs in } w_{i+1} \mbox{ once }, \cr
                j &  \mbox{ otherwise }. \cr
               \end{cases} 
\]
For $j \geq n+1$, define
\[
 r_{w_{i+1}}(j)=\begin{cases}
                j+1, & \mbox{ if } s_{2n-j} \mbox{ occurs in } w_{i+1} \mbox{ once }, \cr
                j &  \mbox{ otherwise }. \cr
               \end{cases} 
\]
\textbf{Case 3:} $\mathcal{O}_q(G) =\mathcal{O}_q(\mbox{Spin}(2n))$.\\ In this case, $w_{i+1}$ is either of the form $s_{i}s_{i-1}\cdots s_k$ where $k \leq n-2$ or 
$s_{i}s_{i+1}\cdots s_{n-2}s_{n-1}s_n$ or $s_{i}s_{i+1}\cdots s_{n-2}s_{n-1}s_ns_{n-2}s_{n-3} \cdots  s_k$. Therefore
for all $1\leq k < n-1$, $s_k$ occurs in  $w_{i+1}$ either  once or  twice or does not occur. Moreover, $s_{n-1}s_n$ can occur either once or does not occur.  For $j \leq n-1$, define
\[
 r_{w_{i+1}}(j)=\begin{cases}
                j-1, & \mbox{ if } s_{j-1} \mbox{ occurs in } w_{i+1} \mbox{ once }, \cr
                j &  \mbox{ otherwise }. \cr
               \end{cases} 
\]
For $j > n+1$, define
\[
 r_{w_{i+1}}(j)=\begin{cases}
                j+1, & \mbox{ if } s_{2n-j} \mbox{ occurs in } w_{i+1} \mbox{ once }, \cr
                j &  \mbox{ otherwise }. \cr
               \end{cases} 
\]
Also, $r_{w_{i+1}}(n)$ and $r_{w_{i+1}}(n+1)$ are $n$ and $n+1$ respectively if $s_{n-1}s_n$ does not occur in $w_{i+1}$ and $n+1$ and $n$ respectively if $s_{n-1}s_n$ occur once in $w_{i+1}$.
Now by a direct verification using diagrams, one can establish the claim. \\
\textbf{Proof of part (2):} Since 
\[
 \{M_n^i, M_n^i+1, \cdots  ,N_n^i\} \subset \{M_n^{i+1}, \cdots  ,N_n^{i+1}\} \subset \cdots \subset \{M_n^{l},  \cdots  ,N_n^{l}\}
\]
the map $r_w: \{M_n^{i}, \cdots  ,N_n^{i}\} \rightarrow \{M_n^{l}, \cdots  ,N_n^{l}\}$ 
is well defined. To verify the claim, first note using diagram representation and possible forms of $w_{i+1}$ given in subsection $2.2$ that  $\pi_{w_{i+1}}^n(u_{j}^{k})=0$
for $k \in   \{M_n^i, \cdots  ,N_n^i\} $ and $  j \notin\{M_n^{i+1}, \cdots N_n^{i+1}\})$. Applying this, we get 
\begin{IEEEeqnarray}{lCl}
 \pi_{w_{i+1}w_{i+2}}^n(u_{r_{w_{i+2}}\circ r_{w_{i+1}}(k)}^k)(e_0\otimes e_0\otimes \cdots \otimes e_0) \nonumber \\
 = \sum_{j=1}^n\pi_{w_{i+1}}^n(u_j^k) (e_0\otimes e_0\otimes \cdots \otimes e_0) \otimes \pi_{w_{i+2}}^n(u_{r_{w_{i+2}}\circ r_{w_{i+1}}(k)}^j)(e_0\otimes e_0\otimes \cdots \otimes e_0) \nonumber \\
 =\sum_{j=M_n^{i+1}}^{N_n^{i+1}}\pi_{w_{i+1}}^n(u_j^k) (e_0\otimes e_0\otimes \cdots \otimes e_0) \otimes \pi_{w_{i+2}}^n(u_{r_{w_{i+2}}\circ r_{w_{i+1}}(k)}^j)(e_0\otimes e_0\otimes \cdots \otimes e_0) \nonumber \\
 = \pi_{w_{i+1}}^n(u_{r_{w_{i+1}}(k)}^k) (e_0\otimes e_0\otimes \cdots \otimes e_0)  \otimes \pi_{w_{i+2}}^n(u_{r_{w_{i+2}}\circ r_{w_{i+1}}(k)}^{r_{w_{i+1}}(k)})(e_0\otimes e_0\otimes \cdots \otimes e_0) \nonumber \\
  \quad \qquad \qquad \qquad \qquad \qquad \qquad   \qquad   \quad( \mbox{by equation } (\ref{m2}))\nonumber \\
 = C(e_0\otimes e_0\otimes \cdots \otimes e_0)  \quad \qquad \quad \qquad   (\mbox{by equation } (\ref{m1}))\nonumber \\\nonumber 
\end{IEEEeqnarray}
Proceeding in this way, one can verify equation (\ref{m1}) for $w=w_{i+1}w_{i+2}\cdots w_{l}$. 
Equation (\ref{m2}) can be verified similarly.
\qed
 \brmrk
 Observe that  the above lemma is not true if $G$ is of type $B$. In particular, from the vertex $n+1$, there is no unique vertex having the property mentioned in Lemma \ref{unique path}.  That's why we excluded 
 type $B$ case.
 \ermrk 

In the following lemma, we define some endomorphisms for each $1\leq i < n$ that act only on $i^{th}$ part $V_{w_i}$ of the simple module $V_w$ associated with the  Weyl word $w$ and  keep other parts unaffected. 
Moreover, on the  $i^{th}$ part, they act 
exactly in the  same way as the variables given in Lemma \ref{sample1} assuming 
that the rank of  the Lie algebra $\mathfrak{g}$ associated with $\mathcal{O}_q(G)$ given there is $i$. 
Here we  are considering two Lie algebras of same type but of different ranks simultaneously. 
So, instead of $\mathfrak{g}$,  we will use  its rank in the notation  $(\!(u_{j,\mathfrak{g}}^i)\!)$ to avoid any confusion. 
  \blmma \label{operator}
 Assume that $\mathcal{O}_q(G)$ is  one of the Hopf $*$-algebras $\mathcal{O}_q(SU(n+1))$, $\mathcal{O}_q(SP(2n))$ or $\mathcal{O}_q(\mbox{Spin}(2n))$.
   For $w \in W_n$, let  $V_w=V_{w_1}\otimes V_{w_2}\otimes \cdots \otimes V_{w_n}$ be the associated 
$\mathcal{O}_q(G)$-module. Then for each $1 \leq i < n$, there exist endomorphisms 
   $T_1^i,T_{2}^i, \cdots T_{N_i}^i  \in \{\pi_w^n(u_{l,n}^k): 1 \leq l,k \leq N_n\}$ such that 
   for all $1\leq j \leq N_i$, one has
  \begin{IEEEeqnarray}{rCl}
   T_j^i(v_1\otimes v_2 \otimes e_0 \otimes \cdots \otimes e_0)&=&C (\pi_{w_1w_2 \cdots w_i}^{i}(u_{j,i}^{N_i})(v_1\otimes v_2) ) \otimes e_0 \otimes \cdots \otimes e_0 \nonumber \\
   &=& Cv_1\otimes  (\pi_{w_i}^{i}(u_{j,i}^{N_i})v_2)  \otimes e_0 \otimes \cdots \otimes e_0 \nonumber 
  \end{IEEEeqnarray}
 where $v_1 \in V_{w_1}\otimes V_{w_2}\otimes \cdots \otimes V_{w_{i-1}}$, $v_2 \in V_{w_i}$ and $C$ is a nonzero constant. 
  \elmma
  \prf
   Define the endomorphism
   \[
    T_j^i= \pi_w^n(u_{r_{w_{i+1}w_{i+2}\cdots w_{n}}(j),n}^{N_n^i})
   \]
for $1 \leq j \leq N_i$.
  Here $r_{w_{i+1}w_{i+2}\cdots w_{n}} : \{M_n^i,M_n^i+1,\dots ,N_n^i\} \rightarrow \{1,2,\dots ,N_n\}$   
  be the function defined in Lemma \ref{unique path}. First note using diagram representation  and possible forms of $w_1\cdots w_i$ given in subsection $2.2$  that 
  $\pi_{w_1\cdots w_i}^n(u_{k}^{N_n^i})=0$ for $k \notin \{M_n^i, \cdots N_n^i\}$. Take 
  $v_1 \in V_{w_1}\otimes V_{w_2}\otimes \cdots \otimes V_{w_{i-1}}$ and  $v_2 \in V_{w_i}$. Then one has  
  \begin{IEEEeqnarray}{lCl}
   \pi_w^n(u_{r_{w_{i+1}w_{i+2}\cdots w_{n}}(j),n}^{N_n^i})(v_1\otimes v_2 \otimes e_0 \otimes \cdots \otimes e_0) \nonumber \\
   =\sum_{k=1}^n\pi_{w_1\cdots w_i}^n(u_{k,n}^{N_n^i})(v_1\otimes v_2) \pi_{w_{i+1}w_{i+2}\cdots w_{n}}^n(u_{r_{w_{i+1}w_{i+2}\cdots w_{n}}(j),n}^{k})
   (e_0 \otimes \cdots \otimes e_0) \nonumber \\
  =\sum_{k=M_n^i}^{N_n^i}\pi_{w_1\cdots w_i}^n(u_{k,n}^{N_n^i})(v_1\otimes v_2) \pi_{w_{i+1}w_{i+2}\cdots w_{n}}^n(u_{r_{w_{i+1}w_{i+2}\cdots w_{n}}(j),n}^{k})(e_0 \otimes \cdots \otimes e_0) \nonumber \\
   =C (\pi_{w_1w_2 \cdots w_i}^{i}(u_{j,i}^{N_i})(v_1\otimes v_2) ) \otimes e_0 \otimes \cdots \otimes e_0  \quad \qquad 
   (\mbox{ by Lemma } \ref{unique path}) \nonumber \\
   = Cv_1\otimes (\pi_{w_i,}^{i}(u_{j,i}^{N_i})(v_2))\otimes e_0 \otimes \cdots \otimes e_0. \nonumber 
  \end{IEEEeqnarray}
 This settles the claim.
\qed\\
The following lemma is an extension of the result of Lemma \ref{sample1}  incorporating all parts of the  Weyl word. 
\blmma \label{poly}
 Assume that $\mathcal{O}_q(G)$ is  one of the Hopf $*$-algebras $\mathcal{O}_q(SU(n+1))$, $\mathcal{O}_q(SP(2n))$ or $\mathcal{O}_q(\mbox{Spin}(2n))$.
For $w \in W_n$, let  $V_w$ be the associated   $\mathcal{O}_q(G)$-module. Then for $1 \leq i \leq n$ and $1 \leq j \leq \ell(w_i)$,  there exist polynomials $p_{j}^{(w,i)}$ with noncommuting variables $\pi_w^n(u_s^r)$'s and
 permutations $\sigma_i$ of $\{1,2,\cdots ,\ell(w_i)\}$   such that 
\begin{IEEEeqnarray}{lCl}
   (p_{1}^{(w,n)})^{r_{\sigma_n(1)}^n} (p_{2}^{(w,n)})^{r_{\sigma_n(2)}^n}\cdots (p_{\ell(w_n)}^{(w,n)})^{r_{\sigma_n(\ell(w_n))}^n}\cdots (p_{1}^{(w,1)})^{r_{\sigma_1(1)}^1}
   \cdots   (p_{\ell(w_i)}^{(w,1)})^{r_{\sigma_1(\ell(w_i))}^1}(e_0\otimes e_0 \otimes \cdots 
   \nonumber \\
 \cdots  \otimes e_0)= C  e_{r_{1}^1}\otimes e_{r_{2}^1} \otimes  \cdots \otimes  e_{r_{\ell(w_1)}^1}  \cdots \otimes e_{r_{1}^n}\otimes
        e_{r_{2}^n} \otimes
 \cdots \otimes e_{r_{\ell(w_n)}^n}  \nonumber 
\end{IEEEeqnarray}
where $C$ is a nonzero constant.
\elmma
\prf For $1 \leq j \leq \ell(w_n)$, let  $p_j^{(w,n)}$ be the polynomial as given in Lemma \ref{sample1} and  $\sigma_n$ be the associated permutation.
To define permutations $\sigma_i$ and polynomials $p_j^{(w,i)}$ for $1\leq i<n$ and $1 \leq j \leq \ell(w_i)$, we view $w_1w_2\cdots w_i$ as an element of the Weyl group of $G$ of rank $i$. So we 
can define polynomial $p_j^{(w_1w_2\cdots w_i,i)}$ and the permutation $\sigma_i$ from Lemma \ref{sample1}. 
Replace the variables $\pi_{w_1w_2\cdots w_i}^i(u_{k,i}^{N_i})$ with $T_k^i$ for $1 \leq k \leq N_i$  in the polynomial $p_j^{(w_1w_2\cdots w_i,i)}$ to define the polynomial 
$p_{j}^{(w,i)}$ for all $1\leq j \leq \ell(w_i)$. 
Now the claim follows from Lemma \ref{sample1}   and Lemma \ref{operator}. 
\qed 
\bthm \label{GKmodule}
 Assume that $\mathcal{O}_q(G)$ is  one of the Hopf $*$-algebras $\mathcal{O}_q(SU(n+1))$, $\mathcal{O}_q(SP(2n))$ or $\mathcal{O}_q(\mbox{Spin}(2n))$.
For $w \in W_n$ and $t \in \bbbt^{n}$, let $V_{w,t}$ be the associated simple unitarizable left  $\mathcal{O}_q(G)$-module. Then one has
\[
 \mbox{GKdim}(V_{w,t})=\ell(w).
\]
\ethm
\prf
Since GKdim$(V_{w,t}) =$ GKdim$(V_{w})$,  we will consider the $\mathcal{O}_q(G)$-module $V_w$ only. Define
\[
 M_0:=\max\{\mbox{degree of } p_j^{(w,i)}: 1\leq i\leq n, 1\leq j\leq \ell(w_i)\}.
\]
 Then by Lemma \ref{poly}, we have
\[
 \mbox{span}\pi_{w}^n(\xi^{M_0k})(e_0\otimes \cdots \otimes e_0) \supset \mbox{span}\{e_{\alpha_1} \otimes e_{\alpha_2} \otimes \cdots e_{\alpha_{\ell(w)}}: \sum_{i=1}^{\ell(w)}\alpha_i=k\}.
\]
Hence dim $\mbox{span }\pi_{w}^n(\xi^{M_0k})(e_0\otimes \cdots \otimes e_0) \geq {k+\ell(w)-1\choose k}$. Therefore
\[
  \mbox{GKdim}(V_w)\geq \varlimsup \frac{\ln \dim(\xi^{M_0k}(e_0\otimes \cdots \otimes e_0))}{\ln M_0k} \geq \varlimsup\frac{\ln {k+\ell(w)-1\choose k}}{\ln M_0k} =\ell(w).
\]
This together with Lemma \ref{upper} proves the claim.
\qed \\
Consider  $V_{w,t}$ as a simple unitarizable left  $\mathcal{O}_q(SO(2n))$-module by restricting the action of $\mathcal{O}_q(\mbox{Spin}(2n))$ to $\mathcal{O}_q(SO(2n))$. 
\bcrlre
For $w \in W_n$ and $t \in \bbbt^{n}$, let $V_{w,t}$ be the associated simple unitarizable left $\mathcal{O}_q(SO(2n))$-module. Then one has
\[
 \mbox{GKdim}(V_{w,t})=\ell(w).
\]
\ecrlre
\prf By Theorem \ref{GKmodule} and the fact that  $\mathcal{O}_q(SO(2n))$ is a subalgebra of $\mathcal{O}_q(\mbox{Spin}(2n))$, it follows that  $\mbox{GKdim}(V_{w,t}) \leq  \ell(w)$. Further observe that 
the variables involved in the polynomials $p_j^{(w,i)}$'s  are entries of the matrix $(\!(u_j^i)\!)$ and hence in $\mathcal{O}_q(SO(2n))$. Therefore Lemma \ref{poly} holds for $\mathcal{O}_q(SO(2n))$-module $V_{w,t}$ 
and the same proof as given in Theorem \ref{GKmodule} will work to show the equality.
\qed

\newsection{Quotient spaces}
Fix a subset $S \subset \Pi$. Let $\mathcal{O}_q(G/K_{S})$ be the quotient Hopf $*$-subalgebra  of 
$\mathcal{O}_q(G)$ (see equation (4.4), \cite{StoDij-1999aa}) . 
If $S$ is the empty set $\phi$, define $W_{\phi}=\{id\}$. 
For a nonempty set $S$, define $W_S$ to be the subgroup of $W_n$ generated by the  simple reflections $s_{\alpha}$ with $\alpha \in S$. 
Let
\[
 W^S:=\{w \in W_n : \ell(s_{\alpha}w)>\ell(w) \quad \forall \alpha \in S\}.
\]
Dijkhuizen and  Stokman (\cite{StoDij-1999aa}) classified  simple unitarizable $\mathcal{O}_q(G/K^{S})$-modules 
in terms of elements of $W^S$. To avoid unnecessary 
complications, let us denote the action of  $\mathcal{O}_q(G/K_S)$ arising by restricting 
 $\pi_w^n$ to $\mathcal{O}_q(G/K_S)$  and the associated module $V_w$
by the same notations $\pi_{w}^n$  and 
$V_w$ respectively. 
\bthm \label{allsimplequotient} (Theorem 3.5, \cite{Sto-2003aa})
The set $\left\{V_{w}; w \in W^S\right\}$ is a complete set of mutually inequivalent simple 
unitarizable left $\mathcal{O}_q(G/K_S)$-module.
\ethm
Let $S_1$ be the empty subset of $\Pi$.
For $2 \leq m \leq n$,   define $S_m$ to be the set  $\{1,2,\cdots ,m-1\}$ if 
$\mathcal{O}_q(G)=\mathcal{O}_q(SU{(n+1)})$ and $\{n-m+2, \cdots ,n\}$ if $\mathcal{O}_q(G)= \mathcal{O}_q(SP(2n))$.
\bthm \label{quomodule}
Let  $\mathcal{O}_q(G)$ be one of the Hopf $*$-algebra $\mathcal{O}_q(SU(n+1))$ or $\mathcal{O}_q(SP(2n))$.
 For $w \in W^{S_{n-m+1}}$,  let $V_{w}$ be the associated simple unitarizable left  $\mathcal{O}_q(G/K_{S_{n-m+1}})$-module. Then one has
\[
 \mbox{GKdim}(V_{w})=\ell(w).
\]
\ethm
\prf  Take  finite sets $M \subset \mathcal{O}_q(G/K_{S_{n-m+1}})$ and 
$F \subset V_{w}$. 
Following the same arguments as given in Lemma \ref{upper}, one can show that 
\[
 \varlimsup\frac{\ln \dim(M^rF)}{\ln r} \leq \ell(w).
\]
Therefore 
\[
 \mbox{GKdim}(V_{w})=sup_{F,M}\varlimsup\frac{\ln \dim(M^rF)}{\ln r} \leq \ell(w).
\]
To show the equality,  observe that 
\begin{enumerate}
 \item 
 for an element $w \in W^{S_{n-m+1}}$  and $1 \leq r \leq n-m$, the $r^{th}$-part
 $w_r=\psi_{r,k_{r}}^{(\epsilon_{r})}(w)$ is identity element of $W_n$. Hence $w$ can be written uniquely as $w= w_{n-m+1} w_{n-m+2}\cdots w_{n}$.
 \item 
 It follows from the definition  that  entries of last 
$m$ rows of $(\!(u_j^i)\!)$ is in the quotient algebra $\mathcal{O}_q(G/K_{S_{n-m+1}})$.
 \item
 The polynomials $p_j^{(w,k)}$ for $1 \leq j \leq \ell(w_k)$ and $n-m+1 \leq k \leq n$ involve variables consisting of entries of last 
$m$ rows of $(\!(u_j^i)\!)$. 
\end{enumerate}
With these facts,  the same arguments  used in Theorem \ref{GKmodule} will prove the claim.                                                                                                              
\qed

\noindent{\sc Partha Sarathi Chakraborty} (\texttt{parthac@imsc.res.in})\\
         {\footnotesize Institute of Mathematical Sciences (HBNI),  CIT Campus,
Taramani,
Chennai, 600113, INDIA}

\noindent{\sc Bipul Saurabh} (\texttt{saurabhbipul2@gmail.com})\\
         {\footnotesize Institute of Mathematical Sciences (HBNI), CIT Campus,
Taramani,
Chennai, 600113, INDIA}

\end{document}